\documentclass{amsart}
\usepackage{amsmath}
\usepackage{amscd}
\usepackage{amssymb}

\newtheorem{theorem}{Theorem}[section]
\newtheorem{theorem/definition}{Theorem/Definition}[section]
\newtheorem{proposition}{Proposition}[section]
\newtheorem{lemma}{Lemma}[section]

\theoremstyle{remark}
\newtheorem{remark}{Remark}[section]
\theoremstyle{definition}

\begin{document}
\title
{On the structure of gradient Yamabe Solitons}

\author{Huai-Dong Cao, Xiaofeng Sun and Yingying Zhang}
\address{Department of Mathematics\\ Lehigh University\\
Bethlehem, PA 18015} \email{huc2@lehigh.edu\,;  xis205@lehigh.edu\,; yiz308@lehigh.edu}

\begin{abstract} We show that every complete nontrivial gradient Yamabe soliton admits a special global warped product structure with a one-dimensional base. 
Based on this, we prove a general classification theorem for complete nontrivial locally conformally flat gradient Yamabe solitons.

\end{abstract}

\maketitle
\date{}

\footnotetext[1]{The work of the first author was partially supported by NSF Grant
DMS-0909581.}

\footnotetext[2]{ The work of the second author was partially supported by NSF Grant DMS-1006696.}

\section{The results}

Self-similar solutions and translating solutions, often called soliton solutions, have emerged in recent years as important objects
in geometric flows since they appear as possible singularity models. Much progress has been made recently in the study of soliton solutions
of the Ricci flow (i.e. Ricci solitons) and the mean curvature flow. In this paper we are interested in geometric structures of Yamabe solitons,
which are soliton solutions to the Yamabe flow. Note that the Yamabe flow has been studied extensively in recent years, see, e.g., 
the very recent survey by Brendle \cite{Bre3} and the references therein. 

\smallskip
A complete Riemannian metric $g=g_{ij}dx^idx^j$ on a smooth manifold $M^n$ is called a {\it gradient  Yamabe soliton} if there exists a smooth
function $f$ such that its Hessian satisfies the equation
$$
\nabla_i\nabla_jf =  (R-\rho) g_{ij},  \eqno(1.1)
$$
where $R$ is the scalar curvature of $g$ and $\rho$ is a constant. For $\rho=0$ the Yamabe soliton is {\it
steady}, for $\rho>0$ it is {\it shrinking} and for $\rho<0$ {\it
expanding}. The function $f$ is called a {\it potential function}
of the gradient Yamabe soliton. When $f$ is constant, we call it a {\sl trivial} Yamabe soliton. It has been known (see \cite{CLN, DD, DS, hsu})  that
every compact Yamabe soliton is of constant scalar curvature, hence trivial since $f$ is harmonic and thus is constant. 

Recently, motivated by the classification of locally conformally flat Ricci solitons and especially \cite{caochen1},  Daskalopoulos and Sesum \cite{DS} initiated
the investigation of conformally flat Yamabe solitons and proved the following very nice classification result:

\begin{theorem} {\bf (Daskalopoulos-Sesum \cite{DS})} \  All complete locally conformally flat gradient Yamabe solitons with positive sectional
curvature $K>0$ are rotationally symmetric.
\end{theorem}
\noindent Moreover, they constructed examples of rotationally symmetric gradient Yamabe solitons on $\mathbb R^n$ with positive
sectional curvature $K>0$.

In this paper, inspired by the above Theorem 1.1 and the recent works of the first author and his collaborators \cite{caochen2, cao+} on Ricci solitons, 
we investigate geometric structures of gradient Yamabe solitons not necessarily locally conformally flat.  It turns out that, by exploring the special nature of the 
Yamabe soliton equation (1.1),  every complete nontrivial gradient Yamabe soliton $(M^n, g, f)$ admits a special global warped product structure with a 1-dimensional base
and the warping function provided by $|\nabla f|$ (see Theorem 1.2).  Based on this special warped product structure, we are able to prove a classification theorem for
locally conformally flat gradient Yamabe solitons without any further assumption on the curvature  (see Corollary 1.5). 
In particular, Theorem 1.1 above is a special case of both Corollary 1.3 and Corollary 1.6(a). 
Our main result is:

\begin{theorem} Let $(M^n, g, f)$ be a nontrivial complete gradient Yamabe soliton satisfying equation (1.1).
Then $|\nabla f|^2$ is constant on regular level surfaces of $f$, and  either

\smallskip
(i)  $f$ has a unique critical point at some point $x_0\in M^n$,  and $(M^n, g, f)$
is rotationally symmetric and equal to the warped product
$$(\mathbb [0,\infty), \, dr^2) \times_{|\nabla f|} (\mathbb S^{n-1},  \, \bar g_{can}),
$$
where $\bar g_{can}$ is the round metric on $\mathbb S^{n-1}$, or

\smallskip
(ii) $f$ has no critical point and $(M^n, g, f)$ is
the warped product
$$
(\mathbb R, \,dr^2) \times\,_{|\nabla f|} \,(N^{n-1}, \, \bar g),
$$ where $(N^{n-1}, \,\bar g)$ is a Riemannian manifold of constant scalar curvature, say $\bar R$.
Moreover, if $(M^n, g, f)$ has nonnegative Ricci curvature $Rc\ge 0$ then $(M^n,g)$ is isometric to the Riemannian product $(\mathbb R, \,dr^2) \times\,(N^{n-1},\bar g)$;
if the scalar curvature $R\geq 0$ on $M^n$, then either  $\bar R>0$,  or  $R=\bar R=0$ and $(M^n,g)$ is isometric to the Riemannian product  $(\mathbb R, \,dr^2) \times\,(N^{n-1},\bar g)$.
\smallskip
\end{theorem}

As an immediate consequence of Theorem 1.2, we have 

\medskip
\noindent {\bf Corollary 1.3.} {\sl Let $(M^n, g, f)$ be a nontrivial complete gradient Yamabe soliton with positive Ricci curvature $Rc>0$, then $f$ has exactly one critical point and $(M^n, g, f)$ 
is rotationally symmetric.}
\medskip

\begin{remark} Shortly after the first version of our paper appeared in the arXiv, G. Catino, C. Mantegazza and L.  Mazzieri \cite{CMM} posted a paper on the 
global structure of conformal gradient solitons with nonnegative Ricci tensor in which they also proved Theorem 1.2 under the assumption of nonnegative 
Ricci curvature (see Theorem 3.2 in \cite{CMM}).  We also remark that, as pointed out in \cite{DS}, steady and shrinking Yamabe solitons have nonnegative 
scalar curvatures.
\end{remark}

In the special case when $(M^n, g, f)$ is locally conformally flat, we can say more about the manifold $(N^{n-1}, \bar g)$ in case (ii) of Theorem 1.2. 

\medskip
\noindent {\bf Theorem 1.4.} {\sl Let $(M^n, g, f)$ be a nontrivial complete gradient Yamabe soliton satisfying equation (1.1). Suppose $f$ has no critical point and is
locally conformally flat, then $(M^n, g, f)$ is the warped product
$$
(\mathbb R, \,dr^2) \times\,_{|\nabla f|} \,(N^{n-1}, \, \bar g_{N}),
$$
where $(N^{n-1}, \,\bar g_{N})$ is a space form (i.e., of constant sectional curvature)}.

\medskip
It is clear that Theorem 1.2 and Theorem 1.4 together implies the following classification of locally conformally flat Yamabe solitons:

\medskip
\noindent {\bf Corollary 1.5.} {\sl Let $(M^n, g, f)$ be a nontrivial complete locally conformally flat gradient Yamabe soliton. Then,
$(M^n, g, f)$ is either

\smallskip
(a)  defined on $\mathbb R^n$, rotationally symmetric, and equal to  the warped product
$$ ([0, \infty),\, dr^2) \times_{|\nabla f|} (\mathbb S^{n-1},  \, \bar g_{can}), \qquad \mbox{or}
$$

\smallskip
(b) the warped product
$$
(\mathbb R, \,dr^2) \times\,_{|\nabla f|} \,(N^{n-1}, \, \bar g_{N}),
$$
where $(N^{n-1}, \,\bar g_{N})$ is a space form.}

\medskip

In particular, we have 

 \medskip
\noindent {\bf Corollary 1.6.} {\sl Let $(M^n, g, f)$ be a nontrivial, non-flat, complete, and locally conformally flat gradient Yamabe soliton.

\smallskip
(a) If $(M^n, g, f)$ has nonnegative Ricci curvature $Rc\ge  0$, then $(M^n, g, f)$ is defined on $\mathbb R^n$ and rotationally symmetric. 

\smallskip
(b) If $(M^n, g, f)$ has nonnegative scalar curvature $R\ge 0$ (as in the steady and shrinking cases), then $(M^n, g, f)$ either is defined on $\mathbb R^n$ 
and rotationally symmetric, or is the warped product cylinder
$$
(\mathbb R,\, dr^2) \times_{|\nabla f|} (\mathbb S^{n-1},  \, \bar g_{can})/\Gamma
$$
for some finite group $\Gamma \subset SO(n)$. }

\smallskip
\begin{remark} Note that Theorem 1.1 also follows from Corollary 1.6(b), since by a well-known theorem of Gromoll and Meyer, $K>0$ implies $M^n$ is 
diffeomorphic to $\mathbb R^n$, hence the latter case in Corollary 1.6(b) cannot happen.
\end{remark}

\begin{remark} As we mentioned before, examples of rotationally symmetric gradient Yamabe solitons on $\mathbb R^n$ with positive
sectional curvature $K>0$ have been constructed by Daskalopoulos and Sesum \cite{DS}. On the other hand, in a forthcoming paper\footnote[3]{The paper is now available on arXiv, 
see reference \cite{He}.}, C. He has shown that any complete gradient steady 
Yamabe soliton on $M^n = \mathbb R\times_{\varphi} N^{n-1}$ is necessarily isometric to the Riemannian product with constant $\varphi$ and $N$ being of zero scalar curvature. Moreover, he showed
the existence of complete gradient Yamabe shrinking soliton metrics on $M^n = \mathbb R \times_{\varphi} N^{n-1}$ with $\rho = 1$, $\bar R>0$ and non-constant warping function $\varphi>0$. 
\end{remark}

\smallskip
\noindent {\bf Acknowledgments}. We would like to thank Chenxu He for very helpful discussions. The first author also likes to thank C. Mantegazza and E. Garcia-Rio 
for bringing the papers \cite{CMM} and \cite{OS} to his attention respectively. 

\section{Warped product structures of complete gradient Yamabe solitons}

We shall follow the notations in \cite{caochen1, caochen2, DS}. Let $(M^n, g_{ij}, f)$ be a complete nontrivial Yamabe soliton, satisfying the Yamabe soliton equation
$$
(R-\rho) g_{ij}=\nabla_i\nabla_jf.
$$
For any regular value $c_0$ of the potential function $f$,  consider the level surface $\Sigma_{c_0}=f^{-1}(c_0)$. Suppose $I$ is an open interval containing $c_0$ such that
$f$ has no critical points
in the open neighborhood $U_{I}=f^{-1} (I)$ of $\Sigma_{c_{0}}$. Then we can express the soliton metric $g$ on
 $U_{I}$ as
$$
ds^2=\frac{1}{|\nabla f|^2} df^2+ {\bar g}_{\Sigma_{c_{0}}},
$$
where
$ {\bar g}_{\Sigma_{c_{0}}}=g_{ab}(f, \theta) d\theta^a d\theta^{b} $
is the induced metric and $\theta=(\theta^2, \cdots, \theta^n)$ is  any local coordinates system on $\Sigma_{c_{0}}$.

On the other hand, as shown in \cite{DS}, we have
$$
\nabla (|\nabla f|^2)=2 \nabla^2 f (\nabla f, \cdot)=2(R-\rho) \nabla f. \eqno(2.1)
$$
Hence,  $|\nabla f|^2$ is constant on any regular level surface $\Sigma_c=f^{-1}(c)\subset U_{I}$, which are all diffeomorphic to $\Sigma_{c_{0}}$.
This allows us to make the change of  variable by setting, up to  an additive constant,
$$ r(x)=\int \frac{df} {|\nabla f|}, \eqno(2.2)
$$
so that
we can further express the metric $g$ on $U_{I}$ as
$$ ds^2=dr^2+g_{ab}(r, \theta) d\theta^a
d\theta^{b}. \eqno(2.3)
$$
Let $\nabla r=\frac{\partial} {\partial r}$, then $ |\nabla r|=1$ and $\nabla f=f'(r) \frac{\partial} {\partial r}$ on $U_{I}$.
Note that  $f'(r)$ does not change sign on $U_{I}$ because  $f$ has no critical points there. Thus, we may assume $I=(\alpha, \beta)$
with  $f'(r)>0$ for $r\in (\alpha, \beta)$.  It is also easy to check that
$$ \nabla_{\frac {\partial} {\partial r}} \frac{\partial} {\partial r}=0, \eqno(2.4)
$$
so integral curves to $\nabla r$ are normal geodesics.

Next, by (2.4) and equation (1.1), it follows that
$$
R-\rho = \nabla^2 f (\frac{\partial} {\partial r}, \frac{\partial} {\partial r}) = f''(r).  \eqno(2.5)
$$
Therefore, we conclude immediately that the scalar curvature $R$ is also constant on $\Sigma_c\subset U_{I}$.
Moreover, the second fundamental form of $\Sigma_c$ is given by
$$
h_{ab}=\frac {\nabla_a\nabla_b f} {|\nabla f|} = \frac {f''(r)} {f'(r)} g_{ab}. \eqno(2.6)
$$
In particular, $\Sigma_c$ is umbilical and its mean curvature is given by

$$
H=(n-1)  \frac {f''(r)} {f'(r)}, \eqno(2.7)
$$
which is again constant along $\Sigma_c$.

Now,  we fix a local coordinates system
$$
(x^1, x^2, \cdots, x^n)=(r, \theta^2, \cdots, \theta^n) \eqno(2.8)
$$
in $U_{I}$, where  $(\theta^2, \cdots,
\theta^n)$ is any local coordinates system on the level surface $\Sigma_{c_{0}}$, and indices $a,b, c, \cdots$ range from $2$ to $n$.
Then, computing in this local coordinates system
we obtain that
\begin{align*}
h_{ab}=- <\partial_r,\nabla_a\partial_b>=-<\partial_r, \Gamma^1_{ab}\partial_r>=-\Gamma^1_{ab}.
\end{align*}
But the Christoffel symbol $\Gamma^1_{ab}$ is given by
\begin{align*}
\Gamma^1_{ab}=\frac{1}{2}g^{11}(-\frac{\partial g_{ab}}{\partial
r})=-\frac 1 2 \frac{\partial g_{ab}}{\partial r}
\end{align*}
Hence, we get
$$
\frac{\partial g_{ab}}{\partial r}=2 \frac {f''(r)} {f'(r)} g_{ab}. \eqno(2.9)
$$
Then, it follows easily from (2.9) that
$$g_{ab} (r, \theta)= (f'(r))^2 g_{ab} (r_0, \theta).   \eqno(2.10) $$ Here the level surface $\{r=r_0\}$ corresponds to $\Sigma_{c_0}$.

Therefore we have arrived at the following

\begin{proposition} Let $(M^n, g_{ij}, f)$ be a complete gradient Yamabe soliton, satisfying the soliton equation (1.1), and let $\Sigma_c=f^{-1}(c)$ be a regular level surface.
Then

\smallskip
(a) $|\nabla f|^2$ is constant on $\Sigma_c$;

\smallskip
(b) the scalar curvature $R$ is constant on $\Sigma_c$;

\smallskip
(c) the second fundamental form of $\Sigma_c$ is given by $h_{ab}=\frac{R-\rho} {|\nabla f|} g_{ab}$;

\smallskip
(d) the mean curvature $H= (n-1) \frac{R-\rho} {|\nabla f|}$ is constant on $\Sigma_c$;

\smallskip
(e) in any open neighborhood $U^{\beta}_{\alpha}=f^{-1}\big ((\alpha, \beta)\big)$ of  $\Sigma_c$ in which $f$ has no critical points, the soliton metric $g$
can be expressed as
$$
 ds^2=dr^2+ (f'(r))^2 {\bar g}_{r_0}. \eqno(2.11)
$$
where $(\theta^2, \cdots, \theta^n)$ is  any local coordinates system on $\Sigma_c$ and ${\bar g}_{r_0}=g_{ab}(r_0, \theta) d\theta^a d\theta^{b} $
is the induced metric on $\Sigma_c=r^{-1}(r_0)$.
\end{proposition}

\begin{remark} Proposition 2.1(a)  was observed first by Sesum and Daskalopoulos \cite{DS}; also Proposition 2.1(b)-(d)
were proved in \cite{DS} under the additional assumption that $(M^n, g, f)$ is locally conformally flat.
\end{remark}

\begin{remark} Our proof of Proposition 2.1 was motivated by arguments in \cite{caochen1, caochen2, cao+} for Ricci solitons. After the first version of our paper appeared in the arXiv, 
we learned that equations similar to Eq. (1.1) had been studied long time ago by various people, see, e.g., \cite{Ta} and the references therein. There are also more recent works, e.g.,  
\cite{OS} and  \cite{CC}.  In particular, Cheeger and Colding \cite{CC} presented beautifully a characterization of warped product structures on a Riemannian manifold $M$ in terms of 
solutions to the more general equation
$$
\nabla_i\nabla_jf =  h \, g_{ij}, 
$$
where $h$ is some smooth function on $M$. 
\end{remark}

Next let us investigate the geometry of the regular level surfaces $\Sigma_c$. To do so, we first need the curvature tensor formula of a
warped product manifold $$(M^n, g)=(I, \,dr^2)\times\varphi \,(N^{n-1}, \, \bar g), \eqno(2.12) $$
where $g=dr^2 +\varphi^2 (r) \bar g$. Fix any local coordinates system $\theta=(\theta^2, \cdots, \theta^n)$ on $N^{n-1}$, and choose
$(x^1, x^2, \cdots, x^n)=(r, \theta^2, \cdots, \theta^n),$ as in (2.8) for the local coordinates system on $M$.  From now on indices $a,b,c,d$ range from $2$ to $n$.
Also curvature tensors with bar are the curvature tensors of $(N, \bar g)$.
Now from either direct computations or \cite{Besse, Oneill}, the Riemann curvature tensor of $(M^n, g)$ is given by
\[
R_{1a1b}=-\varphi\varphi'' \,\bar g_{ab}, \quad R_{1abc}=0, \eqno(2.13)
\]
and
\[
R_{abcd}=\varphi^2 \bar R_{abcd} - (\varphi\varphi')^2\, (\bar g_{ac}\bar g_{bd}-\bar g_{ad} \bar g_{bc}). \eqno(2.14)
\]
The Ricci tensor of $(M^n, g)$ is
$$
R_{11}=-(n-1)\frac {\varphi''} {\varphi}, \quad R_{1a}=0 \quad (2\le a\le n), \eqno(2.15)
$$ and
$$
R_{ab}=\bar R_{ab}-\big [(n-2) (\varphi')^2+\varphi\varphi''\big] \bar g_{ab}  \quad (2\le a,b\le n). \eqno(2.16)
$$
Finally the scalar curvatures of $(M^n, g)$ and $(N^{n-1}, \bar g)$ are related by
$$
R=\varphi^{-2} \bar R- (n-1)(n-2)\left ( \frac{\varphi'}{\varphi}\right )^2-2(n-1)\frac{\varphi''}{\varphi}. \eqno(2.17)
$$
From (2.15) and (2.16), we easily see the following basic facts:

\begin{lemma} (a) The ``radial"  Ricci curvature $R_{11}$ depends only on $r$, hence is constant on level surfaces $\{r\} \times N$;

\medskip
(b) $(N, \bar g)$
is Einstein if and only if the eigenvalues of the Ricci tensor, when restricted to $\{r\} \times N$, are the same and depend only on $r$:
$$R_{22}(r, \theta)=\cdots=R_{nn}(r, \theta) =\mu(r).  \eqno(2.18)$$
\end{lemma}

\begin{remark}
Note that $M^n$ is Einstein, with $Rc=\lambda g$, if and only if $N^{n-1}$ is Einstein with $\bar{Rc}=\bar\lambda \bar g$ and
the warping function $\varphi$ solves the first order ODE
$$
{\varphi'}^2+\frac{\lambda} {n-1}\varphi^2=\frac {\bar\lambda}{n-2}.  \eqno(2.19)
$$
More details can be found in~\cite[9.110]{Besse}.
\end{remark}

Now we are ready to finish the proof Theorem 1.2: \\

\noindent {\sl Proof of Theorem 1.2.}  Let $(M^n, g, f)$ be a complete nontrivial gradient Yamabe soliton.
From Proposition 2.1, we know that $|\nabla f|^2$ is a constant on regular level surfaces of $f$.

Set  $N^{n-1}=f^{-1}(c_{0})$ and $\bar g=\bar g_{r_{0}}$ as in Proposition 2.1 for some regular value $c_0$ of $f$.
Then, since the warping function is $f'(r)$, the warped product formula (2.11) in Proposition 2.1 implies that the potential function $f$ has {at most two}
critical values. Thus, formula (2.11) extends to some maximal interval $I_{max}$,  which is either a finite closed interval  $[\alpha_0, \beta_0]$
with $f'(\alpha_0)=f'(\beta_0)=0$,  or a half-line $[\alpha_0, \infty)$ with $f'(\alpha_0)=0$, or $(-\infty, \infty)$.  However, the first case cannot happen,
for otherwise $M^n$ would be compact, but compact Yamabe solitons are trivial as we mentioned in Section 1.
Thus $f$ has at most one critical point and, after a shift in $r$ if necessary, either $I_{max}=[0, \infty)$, or  $I_{max}=(-\infty, \infty)$.

To see that $(N^{n-1}, \bar g)$ has constant scalar curvature, note that for our Yamabe soliton $(M^n, g, f)$,
we have $R=f''(r)+\rho$ and the warping function is $\varphi (r)=f'(r)$. Thus, from (2.17) we get
$$
 \bar R = (f')^{2} R +(n-1)\big [(n-2) (f'')^2+2f' f'''\big], \eqno(2.20)
$$ which does not depend on $\theta$. Therefore $\bar R$ is a constant.

Also, when $f$ has a unique critical point $x_0$, $r(x)$ is simply the distance function $d(x_0, x)$ from $x_0$. So level surfaces of $f$ are
geodesic spheres centered at $x_0$ which are diffeomorphic to $(n-1)$--sphere $\mathbb S^{n-1}$. In addition, by the smoothness of the metric $g$ at $x_0$ we can conclude
that the induced metric $\bar g$ on $N$ is round (see, e.g., Lemma 9.114 in \cite{Besse}).

Finally, assume we are in the case (ii).  Then, $\varphi=f'>0$ on $(-\infty, \infty)$ and $\varphi'=f''=R-\rho$.  If $(M^n, g,f )$ has nonnegative Ricci curvature  $Rc\ge 0$,
then by (2.15) we know $\varphi''\le 0$, so $\varphi$ is a positive and weakly concave function on $\mathbb R$.  Thus $\varphi$ must be a constant and hence $(M^n, g)$ is isomorphic 
to the Riemannian product $(\mathbb R, dr^2) \times (N, \bar g)$. Now assume $R\geq 0$. Again we prove $\bar R>0$ unless $\bar R=0$ and $(M^n,g)$ is the Riemannian product  $(\mathbb R, \,dr^2) \times\,(N^{n-1},\bar g)$. If $\bar R\leq 0$, by (2.20) we know that $\varphi''\leq 0$. So once again $\varphi$ is a positive weakly concave function on $\mathbb R$ hence a constant function. Therefore, 
again by (2.20), we know that $R=\bar R=0$ and $(M^n,g)$ is the Riemannian product  $(\mathbb R, \,dr^2) \times\,(N^{n-1},\bar g)$. 
\qed

\section{Classification of locally conformally flat Yamabe solitons}

Now let us discuss the classification of locally conformally flat gradient Yamabe solitons and prove Theorem 1.4.

\medskip
\noindent{\sl Proof of Theorem 1.4.}  Let $(M^n, g, f)$ be a nontrivial complete locally conformally flat gradient Yamabe soliton such that $f$ has no critical point.
By Theorem 1.2,  $(M^n, g)$ is a warped product
$$(\mathbb R, \, dr^2) \times_{|\nabla f|} (N^{n-1},  \, \bar g).
$$
Clearly, it remains to prove the following

\medskip
{\bf Claim 1:}  $(N^{n-1}, \bar g)$ is a space form.
\medskip

Indeed, Claim 1 was first proved by Daskalopoulos and Sesum \cite{DS} (see Proposition 2.2 and Lemma 2.4(iii) in \cite{DS}) where they used B. Chow's Li-Yau type
differential Harnack for locally conformally flat
Yamabe flow \cite {Chow} to show that property (2.18) holds and then deduced that $(N, \bar g)$ is a space form. Here we present a simple and direct proof based on the warped
product structure of $(M^n, g, f)$.

All we need is the explicit formula of the Weyl tensor $W$ for an arbitrary warped product manifold (2.12) which can be easily deduced from (2.13)-(2.17):

\[
W_{1a1b}=\frac{\bar R}{(n-1)(n-2)}\bar g_{ab}-\frac{1}{n-2}\bar R_{ab}, \eqno(3.1)
\]

$$
W_{1abc}=0, \eqno(3.2)
$$
and
\[
W_{abcd}=\varphi^2 \overline W_{abcd}. \eqno(3.3)
\]
Here $\overline W$ denotes the Weyl tensor of $(N, \bar g)$.

Now $(M^n, g, f)$ is a warped product, with $\varphi =f'$, and is locally conformally flat, thus  $W=0$. From (3.1) and (3.3), we see
that $(N,\bar g)$ is Einstein and $\overline W=0$. Thus $N$ is a space form.

This proves Claim 1 and completes the proof of Theorem 1.4.

\qed

\begin{remark} By (3.1)-(3.3), it is clear that $(N, \bar g)$ is a space form if and only if $(M, g)$ is locally conformally flat.

\end{remark}

\end{document}